\documentclass{amsart}
\def\proc#1{\medbreak\noindent{\it #1}\hspace{1ex}\ignorespaces}
\def\ep{\noindent{\hfill $\Box$}}
\def\ack{\subsubsection*{Acknowledgement.}}%

\usepackage{amsmath,amsfonts,amsthm,amssymb,graphicx,quotes,setspace,fancyhdr}
\newtheorem{theo}{Theorem}[section]
\newtheorem{lemma}{Lemma}[section]

\newtheorem{cor}{Corollary}[section]
\newtheorem{example}{\bfseries Example}
\def \isnatural {\in\mathbb{N}}
\def \isreal {\in\mathbb{R}}
\def \iscomplex {\in\mathbb{C}}
\def \toinfty {\rightarrow\infty}
\def \noverm {\frac{n}{m}}
\def \half {\frac{1}{2}}
\newcommand{\tef}{transcendental entire function}

\newcommand{\ef}{entire function}

\newcommand\qfor{\quad\text{for }}
\newcommand\spw{spider{'}s web}
\newcommand\spws{spiders{'} webs}
\newcommand\logreg{log-regular}
\newcommand\psireg{$\psi$-regular}

\newcommand\psiregfun{regularity function}
\begin{document}
%
%
%
%
%
\title[Entire functions with spiders' webs]{Entire functions for which the escaping set is a spider's web}
\author{D. J. Sixsmith}
\address{Department of Mathematics and Statistics \\
	 The Open University \\
   Walton Hall\\
   Milton Keynes MK7 6AA\\
   UK}
\email{ds7266@student.open.ac.uk}
%
%
\begin{abstract}
We construct several new classes of {\tef}s, $f$, such that both the escaping set, $I(f)$, and the fast escaping set, $A(f)$, have a structure known as a {\spw}. We show that some of these classes have a degree of stability under changes in the function. We show that new examples of functions for which $I(f)$ and $A(f)$ are {\spws} can be constructed by composition, by differentiation, and by integration of existing examples. We use a property of {\spws} to give new results concerning functions with no unbounded Fatou components.
\end{abstract}

%
%
\maketitle
%
%
\section{Introduction}
Let $f:\mathbb{C}\rightarrow\mathbb{C}$ be a {\tef}, and denote by $f^n, \ n\isnatural,$ the $n$th iterate of $f$. The \itshape Fatou set \normalfont $F(f)$ is defined as the set of points $z \iscomplex$ such that $(f^n)_{n\isnatural}$ forms a normal family in a neighborhood of $z$. The \itshape Julia set \normalfont $J(f)$ is the complement of $F(f)$. An introduction to the properties of these sets can be found in \cite{MR1216719}.

This paper concerns the \itshape escaping set \normalfont $I(f)$ and the \itshape fast escaping set \normalfont $A(f)$, introduced respectively by Eremenko \cite{MR1102727}, and Bergweiler and Hinkkanen \cite{MR1684251}. These sets are defined as follows: $$I(f) = \{z : f^n(z)\rightarrow\infty\text{ as }n\rightarrow\infty\},$$ and $$A(f) = \{z : \text{there exists } L \isnatural \text{ such that } |f^{n+L}(z)| \geq M^n(R, f), \text{ for } n \isnatural\};$$ see \cite{fast} for this form of the definition of $A(f)$. Here, $$M(r, f) = \max_{|z|=r} |f(z)|, \qfor r > 0,$$ $M^n(r, f)$ denotes repeated iteration of $M(r,f)$ with respect to the variable $r$, and $R > 0$ can be taken to be any value such that $M(r, f) > r$ for $r \geq R$. For simplicity, we only write down this restriction on $R$ in formal statements of results -- elsewhere this should be assumed to be true.

Rippon and Stallard \cite{fast} gave a detailed account of many properties of $A(f)$. Their arguments were frequently based on properties of the set $$A_R(f) = \{z : |f^n(z)| \geq M^n(R, f), \text{ for } n \isnatural\}.$$ In particular they showed that $A_R(f), A(f)$ and $I(f)$ can have a structure known as a {\spw}, and that if $A_R(f)$ is a {\spw} then so are $A(f)$ and $I(f)$. They defined a {\spw} as follows:
\proc{Definition.}
A set $E$ is an (infinite) {\spw} if $E$ is connected and there exists a sequence $(G_n)_{n\isnatural}$ of bounded simply connected domains with $G_n \subset G_{n+1},$ for $n\isnatural$, $\partial G_n \subset E,$ for $n\isnatural$, and $\cup_{n\isnatural} G_n = \mathbb{C}$.
\medbreak
Functions for which $A_R(f)$ is a {\spw} have a number of strong dynamical properties. For example, if $A_R(f)$ is a {\spw} then~$f$ has no unbounded Fatou components \cite[Theorem 1.5(b)]{fast} and $A(f)^c$ has uncountably many components \cite[Theorem 1.2]{Osborne1}. Hence, it is desirable to determine those functions for which $A_R(f)$ is a {\spw}, and in \cite[Section 8]{fast} several classes of such functions were given. A further class was given by Mihaljevi\'c-Brandt and Peter \cite{poincare}.

In this paper we give additional classes of examples. First, in Section \ref{Sprelim}, we prove several new results concerning regular growth conditions, which we use in later sections. These results may also be of independent interest.

In Section \ref{Scomp}, we demonstrate a technique for constructing new {\tef}s for which $A_R(f)$ is a {\spw} by taking finite compositions of functions that satisfy a minimum modulus condition and a regularity condition.

In Section \ref{Sstab}, we show that in certain circumstances when $A_R(f)$ is a {\spw}, then so is $A_R(P(f(Q(z)), z))$, where $P, Q$ are polynomials, and so also is $A_R(f+h)$, where the entire function $h$ has smaller growth, in some sense, than $f$. These results allow us to construct a large class of functions for which $A_R(f)$ is a {\spw}. They also show that the property of having an $A_R(f)$ {\spw} can be stable under changes in $f$, unlike many other dynamical properties.

In Section \ref{Sfinite}, we establish a technique for constructing a large class of {\tef}s of finite order for which $A_R(f)$ is a {\spw}, by modifying the power series of a {\tef} of finite order. This technique is a generalisation of the method used to construct some of the examples in \cite{fast}. We show that this class of examples can be extended by differentiation or integration. By combining the results of Sections~\ref{Scomp}, \ref{Sstab} and \ref{Sfinite}, we give an unexpectedly simple function for which $A_R(f)$ is a {\spw}.

In Section \ref{Sinfinite}, we present a technique for constructing new {\tef}s, of infinite order and with large gaps in their power series, for which $A_R(f)$ is a {\spw}.

Finally, in Section \ref{Sbounded}, we relate our results to previous work on classes of {\tef}s which have no unbounded Fatou components. 

%
%
Throughout the paper we use the following three facts about the maximum modulus function $M(r,f)$ of a {\tef} $f$. The first two are well known, and the third is given in \cite[Lemma 2.2]{MR2544754}:
\begin{equation}
\label{M2}
\frac{\log M(r,f)}{\log r} \rightarrow\infty \text{ as } r\rightarrow\infty,
\end{equation}
\begin{equation}
\label{M3}
\text{if } k>1 \text{ then } \frac{M(kr,f)}{M(r,f)} \rightarrow\infty \text{ as } r\rightarrow\infty,
\end{equation}
and there exists $R > 0$ such that
\begin{equation}
\label{M4}
M(r^c,f) \geq M(r,f)^c, \qfor r \geq R, \ c>1.
\end{equation}

We also use the minimum modulus function defined by
\begin{equation}
m(r, f) = \min_{|z|=r} |f(z)|, \qfor r > 0.
\end{equation}

Finally we use the following notation for a disc $$B(z_0, r) = \{ z : |z - z_0| < r \}, \qfor z_0 \iscomplex, r> 0.$$

%
%
\section{New results on regularity}
\label{Sprelim}
In this section we set out conditions which ensure that $A_R(f)$ is a {\spw}. Many of these conditions require some form of regularity of growth. We prove several new results concerning forms of regularity of growth, which enable us to construct examples of functions with an $A_R(f)$ {\spw} later in the paper.

A pair of conditions that are together necessary and sufficient for $A_R(f)$ to be a {\spw} were obtained in \cite[Theorem 8.1]{fast}.
\begin{theo}
\label{webtheorem}
Let $f$ be a {\tef} and let $R > 0$ be such that $M(r, f) > r$ for $r \geq R$. Then $A_R(f)$ is a {\spw} if and only if there exists a sequence $(G_n)_{n \geq 0}$ of bounded simply connected domains such that, for all $n \geq 0$,
\begin{equation} \label{webtheorema} G_n \supset  B(0, M^n(R,f)) \end{equation}
and
\begin{equation}  \label{webtheoremb} G_{n+1} \text { is contained in a bounded component of } \mathbb{C}\backslash f(\partial G_n). \end{equation}
\end{theo}

%
%
This result is very general, and so, in order to construct examples, the following, more readily applicable, sufficient conditions for $A_R(f)$ to be a {\spw} were established in \cite[Corollary 8.3]{fast}.
\begin{lemma}
\label{weblemma}
Let $f$ be a {\tef} and let $R > 0$ be such that $M(r, f) > r$ for $r \geq R$. Then $A_R(f)$ is a {\spw} if, for some $m > 1$, \\
(a) there exists $R_0 > 0$ such that, for all $r \geq R_0$, 
\begin{equation}
\text{there exists } \rho\in (r, r^m) \text{ with } m(\rho,f) \geq M(r, f), \text{ and} \\
\end{equation}
(b) $f$ has regular growth in the sense that there exists a sequence $(r_n)_{n\geq 0}$ with 
\begin{equation}
r_n > M^n(R,f) \text{ and } M(r_n,f) \geq r^m_{n+1}, \qfor n\geq 0.
\end{equation}
\end{lemma}

%
%
We use the following condition, which is stronger than the regularity condition of Lemma~\ref{weblemma}(b), in order to construct a new class of functions with an $A_R(f)$ {\spw}.
\proc{Definition.}
A {\tef} $f$ is \itshape{\psireg} \normalfont if, for each $m>1$, there exist an increasing function $\psi_m$ and $R_m>0$ such that, for all $r\geq R_m$,
\begin{equation}
\label{PRdef}
\psi_m(r) \geq r \ \  \text{ and } \ \ M(\psi_m(r),f) \geq (\psi_m(M(r,f)))^m.
\end{equation}
For given $m>1$ we call $\psi_m$ a \itshape {\psiregfun} \normalfont for $f$.
\medbreak

This condition is slightly stronger than one used in \cite[Theorem 5]{MR2544754} in connection with {\tef}s with no unbounded Fatou components. That version did not require the regularity function to be increasing. However, all the {\psiregfun}s used in \cite{fast,MR2544754} are, in fact, increasing. It was shown in \cite[Section 8]{fast} that if $f$ is {\psireg}, then it satisfies Lemma~\ref{weblemma}(b) for all $m>1$.

%
%
We also use the following condition, which is stronger than {\psireg}ity, in order to construct several classes of functions with an $A_R(f)$ {\spw}.
\proc{Definition.}
Let $c>0$. A {\tef} $f$ is \itshape{\logreg}, with constant $c$\normalfont, if the function $\phi(t) = \log M(e^t,f)$ satisfies 
\begin{equation}
\label{AHdef}
\frac{\phi'(t)}{\phi(t)} \geq \frac{1+c}{t}, \qfor \text{large } t.
\end{equation}
\medbreak
This condition was used by Anderson and Hinkkanen in \cite[Theorem 2]{MR1452790}, also in connection with {\tef}s with no unbounded Fatou components.  The name {\logreg} was suggested by Aimo Hinkkanen in a private communication. The condition was also used in \cite[Section 8]{fast} in order to construct classes of functions with an $A_R(f)$ {\spw}.

In \cite[Section 7]{MR2544754} it was shown that if $f$ is {\logreg} with constant $c$, then, for all $m>1$, $f$ is {\psireg} with {\psiregfun} $\psi_m(r) = r^{m^{1/c}}$; see also Lemma~\ref{AHbasiclemma} below. Hence if $f$ is {\logreg}, then it satisfies Lemma~\ref{weblemma}(b) for all $m>1$.

%
%
We now state three new results concerning {\psireg}ity and {\logreg}ity. The first concerns the composition of {\psireg} functions.
\begin{theo}
\label{TPRcomposition}
Let $f_1, f_2, \ldots, f_k$ be {\tef}s. Suppose that, for all $j\in\{1, 2, \ldots, k\}$, $f_j$~is~{\psireg}, with {\psiregfun} $\psi_m$ for $m>1$. Let $g= f_1 \circ f_2 \circ \cdots \circ f_k$. Then, for any $c>1$, $g$ is {\psireg} with {\psiregfun} $c \psi_m$ for $m>1$.
\end{theo}

In particular it follows that {\psireg}ity is preserved under iteration.
\begin{cor}
If $f$ is a {\psireg} {\tef}, then so is $f^n$ for all $n\isnatural$.
\end{cor}

%
%
The second result relates to the composition of {\ef}s, one of which is {\logreg}.
\begin{theo}
\label{TAHcomposition}
Let $f_1, f_2, \ldots, f_k$ be non-constant entire functions such that, for some $j\in\{1, 2, \ldots, k\}$, $f_j$ is a {\logreg} {\tef}. Let $g=f_1 \circ f_2 \circ \cdots \circ f_k$. Then $g$ is {\logreg}.
\end{theo}

In particular it follows that {\logreg}ity is preserved under iteration.
\begin{cor}
If $f$ is a {\logreg} {\tef}, then so is $f^n$ for all $n\isnatural$.
\end{cor}

Note that Theorem~\ref{TPRcomposition} requires all functions to be {\psireg} {\tef}s, whereas Theorem~\ref{TAHcomposition} requires just one to be a {\logreg} {\tef} and the others only to be entire.

%
%
The third result shows that if $f$ is {\logreg}, then so is any {\tef} with similar growth.
\begin{theo}
\label{TAHcoupled}
Let $f$ and $g$ be {\tef}s. If $f$ is {\logreg} and there exist $a_1, a_2 \geq 1$ and $R_0 > 0$ such that
\begin{equation}
\label{l3M}
M(r^{a_1}, g) \geq M(r, f) \ \ \text{and} \ \ M(r^{a_2}, f) \geq M(r, g), \qfor r \geq R_0,
\end{equation}
then $g$ is {\logreg}.
\end{theo}

We need three preparatory lemmas to prove these results. The first lemma is from \cite{AHregular}, and gives a necessary condition and a sufficient condition for $f$ to be {\logreg}.
\begin{lemma}
\label{AHbasiclemma}
Let $f$ be a {\tef}. \\
(a) If $f$ is {\logreg}, with constant $c$, then there is an $R_0 > 0$ such that, if $k>1$ and $d=k^c$, then 
\begin{equation}
\label{logregeq}
M(r^k, f) \geq M(r, f)^{kd}, \qfor r\geq R_0.
\end{equation}
(b) If (\ref{logregeq}) holds for some $d, k>1$ and $R_0 > 0$, then $f$ is {\logreg}.
\end{lemma}

%
%
The second lemma comes from Wiman-Valiron theory, (see, for example, \cite{MR0385095}), which was first used in connection with the escaping set by Eremenko \cite{MR1102727}. We first need to introduce some terminology. Let $g(z) = \sum_{n=0}^\infty a_n z^n$ be a {\tef}. Define
\begin{equation}
\label{mueq}
\mu(r) = \sup_n |a_n| r^n = |a_N| r^N, \ r>0,
\end{equation}
to be the \itshape maximal term \normalfont of the power series, and call $N = N(r)$ the \itshape central index\normalfont; if (\ref{mueq}) holds for several $N$, we take $N(r)$ to be the largest of these. Note that $N(r)$ is increasing and $N(r) \rightarrow\infty$ as $r \rightarrow\infty$. Wiman-Valiron theory uses $\mu(r)$ to give results about the behaviour of $g$ near points $z(r), \ r > 0$, that satisfy
\begin{equation}
\label{zreq}
|z(r)| = r \text{ and } |g(z(r))| = M(r,g).
\end{equation}

A key result of Wiman-Valiron theory is the following.
\begin{lemma}
\label{wvlemma}
Suppose that $g$ is a {\tef} and $\alpha > \half$. For $r>0$, let $z(r)$ be a point satisfying (\ref{zreq}), and define $$D(r) = B\left(z(r), \frac{r}{(N(r))^\alpha}\right), \quad r>0.$$ Then there exists a set $E \subset (0, \infty)$ with $\int_E 1/t \ dt < \infty$ such that, for $r \notin E$ and $z \in D(r)$,
\begin{equation}
\label{wveq}
g(z) = \left(\frac{z}{z(r)}\right)^{N(r)}g(z(r))(1 + \epsilon),
\end{equation}
where $\epsilon=\epsilon(r, z)\rightarrow 0$ uniformly with respect to $z$ as $r\rightarrow\infty, r\notin E$.
In particular, if $r$ is sufficiently large and $r \notin E$, then
\begin{equation}
\label{covereq}
g(D(r)) \supset \{w : |w| = M(r, g)\}.
\end{equation}
\end{lemma}

We use Lemma~\ref{wvlemma} to prove a result on the behaviour of the maximum modulus of the composite of two entire functions.
\begin{lemma}
\label{AHlem2}
Suppose that $f$ is a non-constant {\ef} and $g$ is a {\tef}. Then, given $\nu>1$, there exist $R_0, R_1>0$ such that 
\begin{equation}
\label{AHlem2eq}
M(\nu r,f \circ g) \geq M(M(r,g),f) \geq M(r,f \circ g), \qfor r \geq R_0,
\end{equation}
and
\begin{equation}
\label{AHlem2eq2}
M(\nu r,g \circ f) \geq M(M(r,f),g) \geq M(r,g \circ f), \qfor r \geq R_1.
\end{equation}
\end{lemma}
\proc{Proof.}
We first prove (\ref{AHlem2eq}). Let $\alpha>\half$, and let $N(r)$, $E$ and $D(r)$ be related to~$g$ as in Lemma~\ref{wvlemma}. Note that $E$ has finite logarithmic measure, and $N(r) \rightarrow\infty$ as $r \rightarrow\infty$. Hence, for sufficiently large $r$, there exists $r' \in (r, \frac{\nu+1}{2} r)\backslash E$, with
\begin{equation}
\label{wvdisc1}
D(r') \subset B(0, \nu r) \text{  and  } g(D(r')) \supset \{w : |w| = M(r', g)\}.
\end{equation}
Let $w$ be such that $|w|=M(r',g)$ and $|f(w)|=M(M(r',g),f)$. Then, by (\ref{wvdisc1}), there is a $z \in D(r')$ with $g(z) = w$. Hence $$|(f \circ g)(z)| = M(M(r', g), f) > M(M(r, g), f).$$ The first part of (\ref{AHlem2eq}) now follows, by (\ref{wvdisc1}). The second part of (\ref{AHlem2eq}) is immediate.

Equation (\ref{AHlem2eq2}) follows in the same way if $f$ is transcendental. Otherwise, suppose that $f$ is a polynomial. Then, for sufficiently large $r$,
\begin{equation}
\label{polydisc}
f(B(0, \nu r)) \supset \{w : |w| = M(r, f)\}.
\end{equation}
Let $w$ be such that $|w|=M(r,f)$ and $|g(w)|=M(M(r,f),g)$. Then, by (\ref{polydisc}), there is a $z \in B\left(0, \nu r\right)$ with $f(z) = w$. Hence $$|(g \circ f)(z)| = M(M(r, f), g).$$ The first part of (\ref{AHlem2eq2}) follows. The second part of (\ref{AHlem2eq2}) is immediate.
\ep\medbreak
In passing, we note a related result discussed by Bergweiler and Hinkkanen \cite[Lemma~1]{MR1684251} that, if we also have $g(0)=0$, then $$M(6r, f \circ g) \geq M(M(r, g), f), \qfor r>0.$$

%
%
Now we are ready to prove Theorems \ref{TPRcomposition}, \ref{TAHcomposition} and \ref{TAHcoupled}.
\proc{Proof of Theorem~\ref{TPRcomposition}.}
Suppose that $m>1$. We note first a general result. Suppose that $f$ is {\psireg} with {\psiregfun} $\psi_m$, and let $\lambda>1$. Then, for sufficiently large $r$, by (\ref{M3}) and (\ref{PRdef}),
\begin{equation}
\label{apsieq}
M(\lambda\psi_m(r),f) \geq \lambda^m M(\psi_m(r),f) \geq \left(\lambda\psi_m(M(r,f))\right)^m.
\end{equation}
Hence $\lambda\psi_m$ is also a {\psiregfun} for $f$.

Now, let $a = c^{1/(k-1)} > 1$. Suppose that $k=2$. Then, for sufficiently large $r$,
\begin{align*}
M(a \psi_m(r), f_1 \circ f_2) &\geq M(M(\psi_m(r),f_2),f_1) \ \ \ &\text{by Lemma } \ref{AHlem2} \\
                            &\geq M((\psi_m(M(r,f_2)))^m,f_1) \ \ \ &\text{by } (\ref{PRdef}) \\
                            &\geq M(\psi_m(M(r,f_2)),f_1)^m \ \ \ &\text{by } (\ref{M4}) \\
                            &\geq (\psi_m(M(M(r,f_2)),f_1))^{m^2} \ \ \ &\text{by } (\ref{PRdef}) \\
                            &\geq (\psi_m(M(r,f_1 \circ f_2)))^{m^2} \ \ \ &\text{ since } \psi_m \text{ is increasing} \\
                            &\geq (a \psi_m(M(r,f_1 \circ f_2)))^{m}. \nonumber
\end{align*}
Hence $g$ is {\psireg} with {\psiregfun} $a \psi_m$. Finally, $c \psi_m = a \psi_m$, since $k=2$. A similar argument with $f_1 \circ f_2$ and $f_3$, both of which are {\psireg} with {\psiregfun} $a \psi_m$, by (\ref{apsieq}), gives the result for $k=3$. The proof follows similarly for larger values of $k$.
\ep\medbreak
%
%
\proc{Proof of Theorem~\ref{TAHcomposition}.}
It is sufficient to prove the result for $k=2$. Suppose then that $k=2$ and that $f_2$ is {\logreg}. By Lemma~\ref{AHbasiclemma}(a) applied to $f_2$, there are $k, d>1$ and $r_1>0$ such that
\begin{equation}
\label{AHeq1}
M(r^k, f_2) \geq M(r, f_2)^{kd}, \qfor r \geq r_1.
\end{equation}
Choose $\nu$ such that $1 < \nu < d$, put $k'=k\nu > 1$ and $d' = \frac{d}{\nu} > 1$. Then, for sufficiently large~$r$,
\begin{align*}
M(r^{k'},f_1 \circ f_2) &\geq M(\nu r^k,f_1 \circ f_2) \\
										    &\geq M(M(r^k,f_2),f_1) &\text{by Lemma } \ref{AHlem2} \\
                        &\geq M(M(r,f_2)^{kd},f_1) &\text{by } (\ref{AHeq1}) \\
               			    &\geq M(M(r,f_2),f_1)^{kd} &\text{by } (\ref{M4}) \\
               		      &\geq M(r,f_1 \circ f_2)^{k'd'} &\text{by choice of } k', d'.
\end{align*}
Thus $f_1 \circ f_2$ is {\logreg} by Lemma~\ref{AHbasiclemma}(b). If $f_1$ is {\logreg} but $f_2$ is not, then the proof that $f_1 \circ f_2$ is {\logreg} is very similar.
\ep\medbreak
%
%
\proc{Proof of Theorem~\ref{TAHcoupled}.}
Suppose that $f$ is {\logreg} with constant $c$, and $a_1, a_2$ are as in (\ref{l3M}). Choose $k>1$ sufficiently large that $k^c > a_1a_2$, put $d~=~k^c, \ k'~=~a_1a_2k>1$, and $d' = \frac{d}{a_1a_2} > 1$. Then, for sufficiently large~$r$,
\begin{align*}
M(r^{k'},g) &= M(r^{a_1a_2k}, g) \\
						&\geq M(r^{a_2k},f) \ \        &\text{by } (\ref{l3M}) \\
            &\geq M(r^{a_2}, f)^{kd} \ \ &\text{by } (\ref{logregeq}) \\
            &\geq M(r, g)^{kd} \ \ &\text{by } (\ref{l3M}) \\
            &=M(r, g)^{k'd'} \ \ &\text{by choice of } k', d'.
\end{align*}
Hence $g$ is {\logreg} by Lemma~\ref{AHbasiclemma}(b).
\ep\medbreak
We now prove several useful corollaries of Theorem~\ref{TAHcoupled}. The first relates to the derivatives and integrals of {\logreg} functions.
\begin{cor}
\label{Cderivative}
Let $f$ be a {\tef}. Then $f$ is {\logreg} if and only if $f'$ is {\logreg}.
\end{cor}
\proc{Proof.}
By integration and (\ref{M4}), for sufficiently large $r$, 
\begin{equation}
\label{der1}
M(r^2, f') \geq r M(r, f') + |f(0)| \geq M(r, f).
\end{equation}
On the other hand, by Cauchy's derivative formula, for sufficiently large $r$,
\begin{equation}
\label{der2}
 M(r^2,f) \geq M(2r, f)/ r\geq M(r, f').
\end{equation}
The result follows by Theorem~\ref{TAHcoupled}, with $a_1 = a_2 = 2$.
\ep\medbreak
The remaining corollaries of Theorem~\ref{TAHcoupled} are used later to give stability results about $A_R(f)$ {\spws}. While they could be combined, they are stated separately for clarity. The first concerns addition of a function to a {\logreg} function.
\begin{cor}
\label{CAHadd}
Let $f$ be a {\logreg} {\tef}, and let $h$ be an {\ef}. Suppose that there exist $a \in (0, 1)$ and $R_0 > 0$ such that
\begin{equation}
\label{heq}
a M(r, f) \geq M(r, h), \qfor r\geq R_0.
\end{equation}
Then $g=f+h$ is {\logreg}.
\end{cor}
\proc{Proof.}
We observe that 
\begin{equation}
\label{fandgeq}
(1 + a)M(r, f) \geq M(r, g) \geq (1-a) M(r, f), \qfor r \geq R_0.
\end{equation}
The result now follows by (\ref{M4}), and Theorem~\ref{TAHcoupled} with $a_1 = a_2 = 2$.
\ep\medbreak
Note that, by (\ref{M2}), if $h$ is a polynomial, then (\ref{heq}) is satisfied for any {\tef} $f$ and any $a\in(0,1)$.

The second corollary concerns a case where {\logreg}ity is preserved under multiplication.
\begin{cor}
\label{CAHbyz}
Let $f$ be a {\logreg} {\tef}. Then $g(z) = z f(z)$ is {\logreg}.
\end{cor}
\proc{Proof.}
By (\ref{M4}) and (\ref{M2}), for sufficiently large $r$, $M(r^2, f) \geq M(r, f)^2 \geq M(r, g).$ Also, for sufficiently large $r$, $M(r, g) \geq M(r, f)$. The result follows, by Theorem~\ref{TAHcoupled} with $a_1 = 1$ and $a_2 = 2$.
\ep\medbreak
Our final corollary is quite general.
\begin{cor}
\label{CAHpoly}
Let $f$ be a {\logreg} {\tef}. Let $P(w, z)$ be a polynomial of degree at least one in $w$, and let $Q(z)$ be a polynomial of degree at least one. Then $g(z) = P(f(Q(z)), z)$ is {\logreg}.
\end{cor}
\proc{Proof.}
Suppose that $$P(f(Q(z)), z) = a f(Q(z))^{N_1} z^{N_2} + h(z) = g_0(z) + h(z),$$ where $N_1$ is the highest power of $w$ in $P(w, z)$ and $N_2$ is the highest power of $z$ corresponding to $f(Q(z))^{N_1}$. By Theorem~\ref{TAHcomposition}, the function $z \mapsto a f(Q(z))^{N_1}$ is {\logreg}. By Corollary~\ref{CAHbyz}, $g_0$ is {\logreg}. Since, by (\ref{M2}), we have $$\half M(r, g_0) \geq M(r, h), \qfor \text{large } r,$$ the result follows by Corollary~\ref{CAHadd}.
\ep\medbreak

%
%
\section{Using composition to give functions for which $A_R(f)$ is a {\spw}}
\label{Scomp}
In this section we demonstrate that $A_R(g)$ is a {\spw} if $g = f_1 \circ f_2 \circ \cdots \circ f_k$, and the {\ef}s $f_j, j\in\{1, 2, \ldots, k\}$, satisfy certain conditions. We need a preparatory lemma before we can state the results. This lemma is a generalisation of Lemma~\ref{weblemma}, in which condition (a) is relaxed slightly. This condition was also used, independently, in \cite{poincare}.
%
%
\begin{lemma}
\label{compositionlemma}
Let $f$ be a {\tef} and let $R > 0$ be such that $M(r, f) > r$ for $r \geq R$. Then $A_R(f)$ is a {\spw} if, for some $m > 1$, \\
(a) there exists $R_0>0$ such that, for all $r\geq R_0$, there is a simply connected domain $G=G(r)$ with 
\begin{equation}
\label{compeq1}
B(0, r) \subset G \subset B(0, r^m) \text{ and } |f(z)| \geq M(r, f), \qfor z \in \partial G,
\end{equation}
and \\
(b) $f$ has regular growth in the sense that there exists a sequence $(r_n)_{n\geq 0}$ with 
\begin{equation}
\label{compeq2}
r_n > M^n(R,f) \text{ and } M(r_n,f) \geq r^m_{n+1}, \qfor n\geq 0.
\end{equation}
\end{lemma}
\proc{Proof.}
Let $m$ and $R_0$ be as in (a), and choose $(r_n)_{n\geq 0}$ satisfying (\ref{compeq2}) with $r_n > R_0$ for $n \geq 0$. For each $n\geq 0$, let $G_n=G(r_n)$.

First, by (\ref{compeq1}) and (\ref{compeq2}),
\begin{equation}
G_n \supset B(0, r_n) \supset B(0, M^n(R,f)), \qfor n\geq 0,
\end{equation}
and so $(G_n)$ satisfies (\ref{webtheorema}).

Second, by (\ref{compeq1}) and (\ref{compeq2}), if $z \in \partial G_n$ then $|f(z)| \geq M(r_n,f) \geq r_{n+1}^m$. Thus $f(G_n)$ contains $B(0, r_{n+1}^m)$, since $f$ maps points of $B(0, M^n(R,f))$ into $B(0, M^{n+1}(R,f)) \subset B(0, r_{n+1}^m)$. Now $G_{n+1}$ is contained in $B(0, r_{n+1}^m)$ and so is contained in a bounded component of $\mathbb{C}\backslash f(\partial G_n)$. Thus $(G_n)$ satisfies (\ref{webtheoremb}). Hence, by Theorem~\ref{webtheorem}, $A_R(f)$ is a {\spw}.
\ep\medbreak
We note that if $P$ is a non-constant polynomial, then $P$ satisfies Lemma~\ref{compositionlemma}(a) for every $m>1$, taking $G(r) = B(0, r^\alpha)$, where $\alpha\in(1,m)$, and a suitable $R_0$.

%
%
We now state the main results of this section. The first relates to the composition of {\psireg} functions, and the second relates to the composition of {\ef}s, one of which is {\logreg}.
\begin{theo}
\label{compthpsi}
Let $f_1, f_2, \ldots, f_k$ be {\tef}s which satisfy Lemma~\ref{compositionlemma}(a) for some $m>1$. Suppose that, for all $j\in\{1, 2, \ldots, k\}$, $f_j$~is~{\psireg}, with {\psiregfun} $\psi_m$ for $m>1$. Let $g = f_1 \circ f_2 \circ \cdots \circ f_k$. Then $A_R(g)$ is a {\spw}, where $R > 0$ is such that $M(r,g) > r$ for $r \geq R$.
\end{theo}
\begin{theo}
\label{compth}
Let $f_1, f_2, \ldots, f_k$ be {\ef}s. Suppose that, for all $j\in\{1,2,\ldots,k\}$, $f_j$ satisfies Lemma~\ref{compositionlemma}(a) for some $m>1$. Suppose also that, for some $j\in\{1, 2, \ldots, k\}$, $f_j$ is a {\logreg} {\tef}. Let $g = f_1 \circ f_2 \circ \cdots \circ f_k$. Then $A_R(g)$ is a {\spw}, where $R > 0$ is such that $M(r,g) > r$ for $r \geq R$.
\end{theo}
%
%
We need one further lemma before we can prove these results. This lemma also concerns the composition of {\ef}s. 
\begin{lemma}
\label{Lcomposition}
Let $f_1, f_2, \ldots, f_k$ be {\ef}s. Suppose that, for all $j\in\{1,2,\ldots,k\}$, $f_j$ satisfies Lemma~\ref{compositionlemma}(a) with $m=m_j>1$. Let $g~=~f_1~\circ~f_2~\circ~\cdots~\circ~f_k$. Then $g$ satisfies Lemma~\ref{compositionlemma}(a) with $m=m_1 m_2 \ldots m_k$.
\end{lemma}
\proc{Proof.}
It is sufficient to prove the result for $k=2$. The result is immediate if $f_1$ and $f_2$ are both polynomials. Otherwise, let $m_1$ and $m_2$ be as given. Consider first the case that $f_2$ is a {\tef}. For sufficiently large $r$, let $G_1=G_1(r)$ be a simply connected domain such that 
\begin{equation}
\label{G1eq1}
B(0, M(r,f_2)) \subset G_1 \subset B(0, M(r,f_2)^{m_1}),
\end{equation}
and
\begin{equation}
\label{G1eq2}
|f_1(z)|  \geq  M(M(r,f_2), f_1), \qfor z\in\partial G_1.
\end{equation}
For sufficiently large $r$, let $G_2=G_2(r)$ be a simply connected domain such that 
\begin{equation}
\label{G2eq1}
B(0, r^{m_1}) \subset G_2 \subset B(0, r^{m_1m_2}),
\end{equation}
and
\begin{equation}
\label{G2eq2}
|f_2(z)|  \geq  M(r^{m_1}, f_2), \qfor z\in\partial G_2.
\end{equation}
If $f_2(z) \in \partial G_1$ then, by (\ref{G1eq1}), $|z| \geq r$, and so there is a component $G_3$ of $f_2^{-1}(G_1)$ which contains $B(0, r)$. Note that $G_3$ is simply connected, and $f_2$ is a proper map of $G_3$ to $G_1$. 

If $z \in \partial G_2$ then, by (\ref{G2eq2}) and (\ref{M4}),
\begin{equation}
\label{maxeq}
|f_2(z)| \geq M(r^{m_1},f_2) \geq M(r, f_2)^{m_1}, \qfor \text{large } r.
\end{equation}
If $z \in \partial G_3$ then, by (\ref{G1eq1}), $|f_2(z)| \leq M(r, f_2)^{m_1}$. Hence, by the maximum principle, if $z \in G_3$ then $|f_2(z)| < M(r, f_2)^{m_1}$. Thus $\partial G_2 \cap G_3 = \emptyset$, by (\ref{maxeq}), and so $B(0, r) \subset G_3 \subset B(0, r^{m_1m_2})$, by (\ref{G2eq1}). Also, if $z \in \partial G_3$ then, by (\ref{G1eq2}), 
\begin{equation}
|(f_1 \circ f_2)(z)| \geq M(M(r,f_2), f_1) \geq M(r, f_1\circ f_2).
\end{equation}
Hence $f_1 \circ f_2$ satisfies Lemma~\ref{compositionlemma}(a), with $m=m_1m_2$.

Secondly, we consider the case where $f_2$ is a polynomial. Choose $m'$ such that $m'>m_1$. For sufficiently large $r$, let $G_1$ and $G_3$ be the domains from the first part of the proof, and let $G_2=B(0,r^{m'})$. Since $f_2$ is a polynomial, for sufficiently large~$r$,
\begin{equation*}
|f_2(z)|  \geq  M(r, f_2)^{m_1}, \qfor z\in\partial G_2.
\end{equation*}
As in the first part of the proof, $\partial G_2 \cap G_3 = \emptyset$, and so $B(0, r) \subset G_3 \subset B(0, r^{m'})$. Also, if $z \in \partial G_3$ then $|(f_1 \circ f_2)(z)| \geq M(r, f_1\circ f_2)$. Hence $f_1 \circ f_2$ satisfies Lemma~\ref{compositionlemma}(a), with $m=m'>m_1$, in particular with $m=m_1m_2$.
\ep\medbreak
In particular it follows from Lemma~\ref{Lcomposition} that the property of satisfying Lemma~\ref{compositionlemma}(a) for some $m>1$ is preserved under iteration.
\begin{cor}
If $f$ is a {\tef} that satisfies Lemma~\ref{compositionlemma}(a) for some $m>1$, then so is $f^n$ for all $n\isnatural$.
\end{cor}
%
%
We are now able to prove Theorems \ref{compthpsi} and \ref{compth}.
\proc{Proof of Theorem~\ref{compthpsi}.}
By Lemma~\ref{Lcomposition}, $g$ satisfies Lemma~\ref{compositionlemma}(a) for some $m>1$. By Theorem~\ref{TPRcomposition}, $g$ is {\psireg}, and so satisfies Lemma~\ref{compositionlemma}(b) for all $m>1$. The result follows by Lemma~\ref{compositionlemma}.
\ep\medbreak
%
%
\proc{Proof of Theorem~\ref{compth}.}
As in the proof of Theorem~\ref{compthpsi}, $g$ satisfies Lemma~\ref{compositionlemma}(a) for some $m>1$. By Theorem~\ref{TAHcomposition}, $g$ is {\logreg}, and so satisfies Lemma~\ref{compositionlemma}(b) for all $m>1$. The result follows by Lemma~\ref{compositionlemma}.
\ep\medbreak
Rippon and Stallard \cite{AHregular} show that there are examples of {\psireg} functions that are not {\logreg}. Some of these examples have order less than $\half$, and so satisfy Lemma~\ref{compositionlemma}(a) for some $m>1$; see Lemma~\ref{orderhalflemma}. This shows that there are situations in which Theorem~\ref{compthpsi} can be applied, but not Theorem~\ref{compth}.

Finally, we note that the conditions of Theorem~\ref{compth} are satisfied by many of the examples in \cite[Section 8]{fast}, and all the examples in this paper (see Sections~\ref{Sfinite} and \ref{Sinfinite}).

%
%
\section{Stability of $A_R(f)$ {\spws}}
\label{Sstab}
Many known dynamical properties of a {\tef} $f$ are unstable under relatively small changes in $f$. For example, the functions $f_1(z) = \exp(-z)$, $f_2(z) = f_1(z) + z + 2\pi i - 1$ and $f_3(z) = f_1(z) +z + 1$ all have very different Fatou sets (see, for example, \cite[Section 4]{MR1216719}). In this section we prove results which show that, in certain circumstances, $A_R(f)$ {\spws} can be very stable. The first result concerns composition with polynomials.
%
%
\begin{theo}
\label{stabilitytheorem}
Suppose that $f$ is a {\logreg} {\tef} which satisfies Lemma~\ref{compositionlemma}(a) for some $m>1$. Let $P(w, z)$ be a polynomial of degree at least one in $w$, and let $Q(z)$ be a polynomial of degree at least one.

Let $g(z) = P(f(Q(z)), z)$. Then $A_R(g)$ is a {\spw}, where $R > 0$ is such that $M(r, g) > r$ for $r \geq R$.
\end{theo}
\proc{Proof.}
By Corollary~\ref{CAHpoly}, $g$ is {\logreg} and so satisfies Lemma~\ref{compositionlemma}(b) for all $m>1$. Hence we need only prove that $g$ satisfies Lemma~\ref{compositionlemma}(a) for some $m>1$.

As in the proof of Corollary~\ref{CAHpoly}, let $$g(z) = a f(Q(z))^{N_1} z^{N_2} + \cdots,$$ where $N_1$ is the highest power of $w$ in $P(w, z)$, and $N_2$ is the highest power of $z$ corresponding to $f(Q(z))^{N_1}$. By Lemma~\ref{Lcomposition}, $f \circ Q$ satisfies Lemma~\ref{compositionlemma}(a). Hence, there is an $m>1$ such that, for sufficiently large $r$, there is a simply connected domain $G=G(r)$ with $B(0, r^m) \subset G \subset B(0, r^{m^2})$ and
\begin{equation}
\label{Geq1}
|f(Q(z))| \geq M(r^{m}, f \circ Q), \qfor z\in\partial G.
\end{equation}
Hence, when $z \in \partial G$, for sufficiently large $r$,
\begin{align*}
|g(z)| &\geq \half|a| M(r^m, f \circ Q)^{N_1} r^{N_2}  \ \ \ &\text{by (\ref{Geq1}) and (\ref{M2})} \\
       &\geq 2 |a| M(r, f \circ Q)^{N_1} r^{N_2}       \ \ \ &\text{by (\ref{M4})} \\
       &\geq M(r, g)                                   \ \ \ &\text{by (\ref{M2})}.
\end{align*}
Thus $g$ satisfies Lemma~\ref{compositionlemma}(a) with $m$ replaced by $m^{2}$, so the proof is complete.
\ep\medbreak
The second result concerns addition of an {\ef} to a {\tef} with an $A_R(f)$ {\spw}.
%
%
\begin{theo}
\label{stabilitytheorem2}
Suppose that $f$ is a {\logreg} {\tef} which satisfies Lemma~\ref{compositionlemma}(a) for some $m>1$, and that $h$ is an {\ef}. Suppose also that there exist $a \in (0, 1)$ and $R_0 > 0$ such that
\begin{equation}
\label{anotherheq2}
a M(r, f) \geq M(r^m, h), \qfor r \geq R_0.
\end{equation}
Let $g = f + h$. Then $A_R(g)$ is a {\spw}, where $R > 0$ is such that $M(r, g) > r$ for $r \geq R$.
\end{theo}
\proc{Proof.}
First we note that, for sufficiently large $r$, $a M(r, f) \geq M(r, h)$. Hence, by Corollary~\ref{CAHadd}, $g$ is {\logreg} and so satisfies Lemma~\ref{compositionlemma}(b) for all $m>1$. Thus, by Lemma~\ref{compositionlemma}, it remains to prove that $g$ satisfies Lemma~\ref{compositionlemma}(a) for some $m>1$.

By hypothesis, for sufficiently large $r$, there is a simply connected domain $G=G(r)$ with $B(0, r^m) \subset G \subset B(0, r^{m^2})$ and
\begin{equation}
\label{Geq}
|f(z)| \geq M(r^m, f), \qfor z\in\partial G.
\end{equation}
Thus, when $z \in \partial G$, for sufficiently large $r$,
\begin{align*}
|g(z)| &\geq |f(z)| - |h(z)|          &\text{ } \\
       &\geq (1-a)M(r^m, f)          	&\text{by (\ref{anotherheq2}), and (\ref{Geq})} \\
       &\geq (1+a) M(r, f)            &\text{by (\ref{M4})} \\
       &\geq M(r, g)                  &\text{ }.
\end{align*}
Hence $g$ satisfies Lemma~\ref{compositionlemma}(a) with $m$ replaced by $m^2$, so the proof is complete.
\ep\medbreak
\proc{Remark.}
Using the same method of proof it can be shown that in Theorem~\ref{stabilitytheorem2} the function $h$ can also be of the form $h(z) = f(z)/(z-c)$, where $f(c)=0$.
\medbreak

Finally, we note that the conditions on $f$ in Theorems \ref{stabilitytheorem} and \ref{stabilitytheorem2} are satisfied by many of the examples in \cite[Section 8]{fast}, and all the examples in this paper. It can be shown that these conditions are also satisfied by the functions in \cite{poincare}. So we can produce new functions for which $A_R(f)$ is a {\spw} by taking these known examples and applying Theorems \ref{stabilitytheorem} and \ref{stabilitytheorem2}.

%
%
\section{Functions of finite order for which $A_R(f)$ is a {\spw}}
\label{Sfinite}
In this section we develop a technique which enables us to take a {\tef} of finite order, modify its power series, and produce a class of {\tef}s of finite order for which $A_R(f)$ is a {\spw}. From the exponential function we obtain a class of such functions (Example~\ref{fEx1}) which contains the function
\begin{equation}
\label{spwexample}
f(z) = \half(\cos z^{\frac{1}{4}} + \cosh z^{\frac{1}{4}}) = \sum_{n=0}^\infty \frac{z^{n}}{(4n)!}
\end{equation}
given in \cite[Section 8]{fast}, together with the related functions
\begin{equation}
f(z) = \sum_{n=0}^\infty \frac{z^{pn}}{(qn)!}, \ \ p, q \isnatural, \ \ p/q < \half,
\end{equation}
suggested by Halburd and also mentioned in \cite[Section 8]{fast}. We obtain another class (Example~\ref{fEx3}) from the error function (see \cite[p.297]{standards})
\begin{equation}
\label{erfeq}
\text{erf}{(z)} = \frac{2}{\sqrt{\pi}}\sum_{n=0}^\infty \frac{(-1)^n}{n!(2n+1)}z^{2n+1}.
\end{equation}

We define the order $\rho(f)$ and lower order $\lambda(f)$ of a {\tef}~$f$ by
\begin{equation}
\label{orderdef}
\rho(f) = \limsup_{r\toinfty} \frac{\log\log M(r,f)}{\log r} \ \ \text{ and } \ \ \lambda(f) = \liminf_{r\toinfty} \frac{\log\log M(r,f)}{\log r}.
\end{equation}

We note from, for example, \cite{MR0281918} that if $f(z) = \sum_{n=0}^\infty a_n z^n$, then 
\begin{equation}
\label{order}
\rho(f) = \limsup_{n\toinfty} \frac{n\log n}{\log|a_n|^{-1}}
\end{equation}
and
\begin{equation}
\label{lowerorder}
\lambda(f) = \max_{(n_p)}\liminf_{p\toinfty} \frac{n_p\log n_{p-1}}{\log|a_{n_p}|^{-1}}.
\end{equation}

%
%
We use the following three lemmas, all discussed in \cite[Corollary 8.3 and the following remarks]{fast}. The first is from \cite[p.205]{MR2458805}, and gives a sufficient condition for a {\tef} to be {\logreg}.
\begin{lemma}
\label{logreglemma}
If $f$ is a {\tef} of finite order and positive lower order, then $f$ is {\logreg}.
\end{lemma}
The second is from, for example, \cite[Satz 1]{MR0097532}.
\begin{lemma}
\label{orderhalflemma}
If $f$ is a {\tef} of order less than $\half$, then $f$ satisfies Lemma~\ref{weblemma}(a) for some $m>1$.
\end{lemma}
The third follows from Lemma \ref{logreglemma}, Lemma \ref{orderhalflemma} and Lemma \ref{weblemma}. 
\begin{lemma}
\label{isaspwlemma}
If $f$ is a {\tef} of order less than $\half$ and positive lower order, then $A_R(f)$ is a {\spw}, where $R > 0$ is such that $M(r, f) > r$ for $r \geq R$.
\end{lemma}

We use the following operator to produce classes of functions which satisfy the conditions of Lemma~\ref{isaspwlemma}.
\proc{Definition.}
For $n, m \isnatural$, let $T_{\noverm}$ be defined by
\begin{equation}
T_{\noverm}(f(z)) = \frac{1}{m}\sum_{k=1}^m f(e^{\frac{2\pi ik}{m}} z^{\noverm}),
\end{equation}
where $f$ is an {\ef}, and we choose a consistent branch of the $m$th root for each term in the sum.
\medbreak
If $f$ is a {\tef}, then the $T_{\noverm}$ operator extracts from the power series of $f$ only those terms with exponents which are multiples of~$m$, and these exponents are multiplied by $n/m$ (see (\ref{Teq}) below). For example, if $f(z) = e^z$, then $$T_{\frac{2}{3}}(f(z)) = 1 + \frac{z^2}{3!} + \frac{z^4}{6!} + \cdots.$$ We note in passing that the $T_{\noverm}$ operator has some appealing properties; for example, $T_{\frac{1}{m}} \circ T_{\frac{1}{n}} = T_{\frac{1}{nm}}$ and also $T_{\noverm}(f(z^m)) = f(z^n)$. 

The following result concerns a key property of this operator, namely its effect on the order of a function.
\begin{theo}
\label{finitetheorem}
If $f$ is a {\tef} of order $\rho(f)$ and~$n,m~\in~\mathbb{N}$, then $T_{\noverm}(f)$ is a well-defined entire function of order at most $\noverm\rho(f)$.
\end{theo}
\proc{Proof.}
First, we consider the action of $T_{\noverm}$ on the power series $f(z) = \sum_{l=0}^\infty a_l z^l$. Since we have a consistent choice of the $m$th root, the sum of the complex roots of unity is zero, and with $p=l/m$, we obtain
\begin{equation}
\label{Teq}
T_{\noverm}(f(z)) = \frac{1}{m}\sum_{k=1}^m \sum_{l=0}^\infty a_l e^{\frac{2\pi ikl}{m}} z^{\frac{ln}{m}} 
                  = \sum_{l=0}^\infty a_l z^{\frac{ln}{m}} \sum_{k=1}^m \frac{1}{m} e^{\frac{2\pi ikl}{m}} 
                  = \sum_{p=0}^\infty a_{pm} z^{pn}.
\end{equation}
Hence the value of $T_{\noverm}(f)$ is independent of the choice of the $m$th root, and this power series has infinite radius of convergence.

We deduce from (\ref{order}), with $k = pm$, that
\begin{align*}
\rho(T_{\noverm}(f)) &= \limsup_{p\toinfty} \frac{pn\log pn}{\log|a_{pm}|^{-1}} \\
                     &\leq \limsup_{k\toinfty}\frac{(kn/m)\log (kn/m)}{\log|a_k|^{-1}} \\
                     &= \noverm \limsup_{k\toinfty} \frac{k\log k}{\log|a_k|^{-1}} \\
                     &= \noverm \rho(f),
\end{align*}
as required.
\ep\medbreak
We now seek to use this operator, together with Lemma~\ref{isaspwlemma}, to generate {\tef}s for which $A_R(f)$ is a {\spw}. It is possible for the function $T_{\noverm}(f)$ to be simply a polynomial when $f$ is a {\tef}. For example, if $f(z) = z\exp(z^2)$ then $T_{\frac{1}{2}}(f(z)) = 0$, because the power series of $f$ has only odd powers of $z$ which are eliminated by the $T_{\frac{1}{2}}$ operator.

Even if $T_{\noverm}(f)$ is transcendental, $T_{\noverm}(f)$ may not have positive lower order when $f$ does. For example, if $g$ is a {\tef} of  order less than $1$ and lower order zero, then $f(z) = g(z^2) + z\exp(z^2)$ has both order and lower order $2$, but $T_{\frac{1}{2}}(f(z))=T_{\frac{1}{2}}(g(z^2))=g(z)$ has order less than $1$ and lower order zero, reasoning as in the previous paragraph.

The following lemma gives two sufficient conditions for $T_{\noverm}(f)$ to have positive lower order. 
\begin{lemma}
\label{finiteorderlemma}
Let $f(z) = \sum_{p=0}^\infty a_p z^p$ be a {\tef}, and let $n,m \isnatural$. \\
(a) If
\begin{equation}
\label{lowerthing}
\liminf_{p\toinfty} \frac{p\log p}{\log|a_{pm}|^{-1}} > 0,
\end{equation}
then $T_{\noverm}(f)$ has positive lower order. \\
(b) If $T_{\noverm}(f)$ has positive lower order, and $g(z) = \sum_{p=0}^\infty b_p z^p$ is a {\tef} with  $|b_p| \geq |a_p|$ for $p$ sufficiently large, then $T_{\noverm}(g)$ has positive lower order.
\end{lemma}
\proc{Proof.}
For part (a) we note, by (\ref{Teq}) and with $n_p = np$ in (\ref{lowerorder}), that
\begin{align*}
\lambda(T_{\noverm}(f)) \geq \liminf_{p\toinfty} \frac{np\log (n(p-1))}{\log|a_{pm}|^{-1}}
											  = n \liminf_{p\toinfty}\frac{p\log p}{\log|a_{pm}|^{-1}} > 0.
\end{align*}
Part (b) follows immediately from (\ref{lowerorder}).
\ep\medbreak

%
%
We now give some explicit examples of classes of functions for which $A_R(f)$ is a {\spw}. The first example includes (\ref{spwexample}) as a special case.
\begin{example}
\label{fEx1}
Let $f = T_{\noverm}(g)$, where $g(z) = \exp(z)$ and where $m>2n$. Then $A_R(f)$ is a {\spw}, where $R > 0$ is such that $M(r,f) > r$ for $r \geq R$.
\end{example}
\proc{Proof.}
The exponential function has order $1$, and satisfies (\ref{lowerthing}) for all $m>1$. Thus $f$ has order less than $\half$ by Theorem~\ref{finitetheorem}, and the result follows by Lemma~\ref{finiteorderlemma}(a) and Lemma~\ref{isaspwlemma}.
\ep\medbreak
The second example illustrates the use of both parts of Lemma~\ref{finiteorderlemma}. This result can also be justified by the results of Sections \ref{Scomp} and \ref{Sstab}.
\begin{example}
\label{fEx2}
Let $f = T_{\noverm}(g)$, where $g(z) = z\exp{(z^2)}+\exp(z)$, $m>4n$ and $m$ is odd. Then $A_R(f)$ is a {\spw}, where $R > 0$ is such that $M(r,f) > r$ for $r \geq R$.
\end{example}
\proc{Proof.}
The function $z \mapsto z\exp{(z^2)}$ has order $2$, and satisfies (\ref{lowerthing}) when $m$ is odd. Thus $f$ has order less than $\half$ by Theorem~\ref{finitetheorem}, and the result follows by Lemma~\ref{finiteorderlemma} and Lemma~\ref{isaspwlemma}.
\ep\medbreak
The technique of this section can be applied any {\tef} of finite order, provided its power series satisfies (\ref{lowerthing}) for some $m\isnatural$. We illustrate this with the error function.
\begin{example}
\label{fEx3}
Let $f = T_{\noverm}(g)$, where $g(z) = \text{erf}(z)$, $m>4n$ and $m$ is odd. Then $A_R(f)$ is a {\spw}, where $R > 0$ is such that $M(r,f) > r$ for $r \geq R$.
\end{example}
\proc{Proof.}
By (\ref{erfeq}) and (\ref{order}), $g$ has order $2$, and satisfies (\ref{lowerthing}) when $m$ is odd. Thus $f$ has order less than $\half$ by Theorem~\ref{finitetheorem}, and the result follows by Lemma~\ref{finiteorderlemma}(a) and Lemma~\ref{isaspwlemma}.
\ep\medbreak
Our final example combines earlier results to give an unexpectedly simple function with an $A_R(f)$ {\spw}.
\begin{example}
\label{fEx4}
Let $f(z) = \cos z + \cosh z$. Then $A_R(f)$ is a {\spw}, where $R > 0$ is such that $M(r,f) > r$ for $r \geq R$.
\end{example}
\proc{Proof.}
This follows from Theorem~\ref{compth} and the function in (\ref{spwexample}).
\ep\medbreak
%
%
Our goal in this section has been to produce a class of {\logreg} {\tef}s of order less than $\half$, which, by Lemmas~\ref{orderhalflemma} and~\ref{weblemma}, have an $A_R(f)$ {\spw}. Finally, we show that this class can be extended by differentiation or integration, thus giving a further method of constructing $A_R(f)$ {\spws}.
\begin{theo}
\label{Tdiff}
Let $f$ be a {\logreg} {\tef} of order less than $\half$, and let $g$ be the derivative of $f$ or an integral of $f$. Then $A_R(g)$ is a {\spw}, where $R > 0$ is such that $M(r,g) > r$ for $r \geq R$.
\end{theo}
\proc{Proof.}
Since $g$ has the same order as $f$, and is {\logreg} by Corollary~\ref{Cderivative}, the result follows by Lemmas~\ref{orderhalflemma} and~\ref{weblemma}.
\ep\medbreak

%
%
\section{A function of infinite order with gaps for which $A_R(f)$ is a {\spw}}
\label{Sinfinite}
We recall that a transcendental entire function $f$ has Fabry gaps if $$f(z)=\sum_{k=1}^\infty~a_k z^{n_k}$$ and $n_k/k \rightarrow\infty$ as $k \rightarrow\infty$. By a result of Fuchs \cite{MR0159933}, an {\ef} $f$ of finite order with Fabry gaps satisfies Lemma~\ref{weblemma}(a) for $m>1$. This fact was used by Wang in \cite[Theorem 1]{MR1808648} to describe a class of entire functions with no unbounded Fatou components. Thus if $f$ is also {\logreg} then, by Lemma~\ref{weblemma}, $A_R(f)$ is a {\spw}. (As noted earlier, a {\logreg} {\tef} satisfies Lemma \ref{weblemma}(b) for all $m>1$.) This fact was pointed out by Rippon and Stallard \cite[Theorem 1.9(d)]{fast}, who gave an example of such a function \cite[Example 1]{fast}, shown to be {\logreg} by using Lemma \ref{logreglemma}.

It was also pointed out in \cite{MR1808648} and in \cite[Section 8]{fast} that, by a result of Hayman \cite{MR0306497}, Lemma~\ref{weblemma}(a) holds in the case of certain functions of infinite order with gaps:
\begin{lemma}
\label{lemmahayman}
Let $f(z) = \sum_{k=1}^\infty a_k z^{n_k}$ be a {\tef} where, for some $\alpha>2$,
\begin{equation}
\label{lemmahaymaneq}
n_k > k \log k(\log\log k)^\alpha, \qfor \text{large } k.
\end{equation}
Then $f$ satisfies Lemma~\ref{weblemma}(a) for $m>1$.
\end{lemma}

Wang \cite[Theorem 2]{MR1808648} used this result to show that if $f$ satisfies (\ref{lemmahaymaneq}) and has a property equivalent to {\logreg}ity, then $f$ has no unbounded Fatou components.

Suppose that $g$ is a {\tef} of infinite order generated by omitting terms from the power series of another {\tef}, $f$ say, and $g$ satisfies (\ref{lemmahaymaneq}). If $g$ is also {\logreg}, then $A_R(g)$ is a {\spw}, by Lemma~\ref{weblemma}. If $f$ has infinite order, then it does not seem straightforward to check that such a function $g$ is {\logreg}. In this section we demonstrate a method for achieving this, and then give an explicit example of such a function.

%
%
We start with a general result.
\begin{theo}
\label{infinitetheorem}
Suppose that $f(z) = \sum_{n=0}^\infty a_n z^n$ is a {\logreg} {\tef} and there exists $N_0\isnatural$ such that
\begin{equation}
0 < a_{n+1} \leq a_n, \qfor n\geq N_0.
\end{equation}
Suppose also that $g$ is a {\tef} with
\begin{equation}
g(z) = \sum_{k=1}^\infty a_{n_k} z^{n_k},
\end{equation}
where, for some $M>1$ and $\alpha>2$,
\begin{equation}
\label{nkdecrease}
1 < \frac{n_{k+1}}{n_k} < M, \qfor \text{large } k,
\end{equation}
and
\begin{equation}
\label{gaps}
n_k > k \log k(\log\log k)^\alpha, \qfor \text{large } k.
\end{equation}
Then $g$ is {\logreg} and $A_R(g)$ is a {\spw}, where $R > 0$ is such that $M(r, g) > r$ for $r \geq R$.
\end{theo}
\proc{Proof.}
By Lemma~\ref{lemmahayman}, $g$ satisfies Lemma~\ref{weblemma}(a) for $m>1$. To complete the proof, we use Theorem~\ref{TAHcoupled} to show that $g$ is {\logreg}.

Without loss of generality, by adding a polynomial, we can assume by Corollary~\ref{CAHpoly} that $N_0 = 0$ and (\ref{nkdecrease}) holds for $k\geq 1$. Because $a_n > 0$ for $n\geq 0$, $$M(r, f) = f(r) > g(r) = M(r, g), \qfor r>0.$$ Thus it remains to show that there exist $a>1$ and $R_0>0$ such that 
\begin{equation}
\label{gandfeq}
M(r^a, g) \geq M(r, f), \qfor r\geq R_0.
\end{equation}

Choose $a'>1$ and $K>1$ sufficiently large such that
\begin{equation} 
\label{stineq}
\frac{n_{k+1}}{n_{k}} < \half(1+a') < a', \ \text{ and } \ K^{n_{k}} > n_{k+1} - n_{k},  \qfor k\geq 1.
\end{equation} 
Now let $\mu = \half(a'-1) > 0$, and define
\begin{align}
\label{stAdef1}
&M(r^{a'}, g) = \sum_{k=1}^{\infty} A_k, \ \ \ A_k = a_{n_k} r^{a'n_k}, \\
\label{stAdef2}
&M(r, f)      = \sum_{n=0}^{a_{n_1}-1} a_n r^n + \sum_{k=1}^{\infty} B_k, \ \ \ B_k = a_{n_{k}} r^{n_{k}} + \cdots + a_{n_{k+1}-1} r^{n_{k+1}-1}.
\end{align}
Because the $a_n$ are decreasing,
\begin{equation}
B_k < \left(n_{k+1}-n_{k}\right)a_{n_{k}}r^{n_{k+1}}, \qfor r>1 \text{ and } k\geq 1.
\end{equation}
Thus, if $r > \max\{1, K^{\frac{1}{\mu}}\}$, then, by (\ref{stAdef1}) and (\ref{stineq}),
\begin{equation}
B_k < \left(n_{k+1}-n_k\right)r^{n_{k+1}-a'n_k}A_k < K^{n_k} r^{-n_k\mu} A_k < A_k,  \qfor k\geq 1.
\end{equation}
Thus, by (\ref{stAdef1}) and (\ref{stAdef2}),
\begin{equation}
\label{zzz}
M(r^{a'}, g) > M(r, f) - \sum_{n=0}^{a_{n_1}-1} a_n r^n, \qfor r > \max\{1, K^{\frac{1}{\mu}}\}.
\end{equation}
Finally, for any $a > a'$ we can choose $r$ sufficiently large such that
\begin{align}
M(r^{a}, g) &\geq 2M(r^{a'}, g) \ \ \ &\text{ by (\ref{M3})} \\
            &> 2M(r, f) - 2\sum_{n=0}^{a_{n_1}-1} a_n r^n \ \ \ &\text{by (\ref{zzz})} \\
            &\geq M(r, f) \ \ \ &\text{ by (\ref{M2})}.
\end{align}
This proves (\ref{gandfeq}) as required.
\ep\medbreak
In the rest of this section we construct an explicit example of a {\tef} $f$ of infinite order, defined by a gap series, for which $A_R(f)$ is a {\spw}. First we need a simple result about functions of infinite order.
\begin{lemma}
\label{inforderlemma}
Let $f$ and $g$ be {\tef}s, and suppose that $f$ has infinite order. If there exist $a, R_0>0$ such that $$M(r^a, g) \geq M(r, f), \qfor r\geq R_0,$$ then $g$ has infinite order.
\end{lemma}
\proc{Proof.}
By (\ref{orderdef}),
\begin{align*}
\rho(g) = \limsup_{r\toinfty} \frac{\log\log M(r^a,g)}{\log r^a} 
        \geq \frac{1}{a} \limsup_{r\toinfty} \frac{\log\log M(r,f)}{\log r} 
        = \frac{1}{a} \rho(f),
\end{align*}
and the result follows.
\ep\medbreak
%
%
The next lemma is needed in the construction of our example.
\begin{lemma}
\label{decreasinglemma}
Let $g(z) = \sum_{n=0}^\infty a_n z^n$ be a {\tef}, with $a_n \isreal$ for $n\geq 0$, $a_1 \leq 1$, and
\begin{equation}
\label{gcons}
0 < (n+1) a_{n+1} \leq n a_n, \qfor n \geq 1.
\end{equation}
Then $f(z) = \exp(g(z))$ has power series $f(z) = \sum_{n=0}^\infty b_n z^n$, where
\begin{equation}
\label{dec1}
0 < b_{n+1} \leq b_n, \qfor n \geq 1.
\end{equation}
\end{lemma}
\proc{Proof.}
Clearly $b_n>0$ for $n\geq 0$. Since $f'(z) = g'(z) f(z)$ we have
\begin{equation}
\sum_{n=0}^\infty (n+1) b_{n+1} z^{n} = \sum_{k=0}^\infty (k+1) a_{k+1} z^k \sum_{l=0}^\infty b_l z^l.
\end{equation}
Equating powers of $z$ gives
\begin{equation}
\label{baeq1}
(n+1) b_{n+1} = \sum_{l=0}^n (n+1-l) a_{n+1-l} b_l, \qfor n \geq 0.
\end{equation}
Hence, for $n\geq 1$,
\begin{align}
\label{baeq2}
(n+1) b_{n+1} &=     \sum_{l=0}^{n-1} (n+1-l) a_{n+1-l} b_l  + a_1 b_n &\text{ } \\
              &\leq  \sum_{l=0}^{n-1} (n-l) a_{n-l} b_l  + b_n,  &\text{by (\ref{gcons}) and as } a_1 \leq 1 \\
              &=     n b_n + b_n, &\text{by } (\ref{baeq1}),
\end{align}
which proves that (\ref{dec1}) holds.
\ep\medbreak
%
%
Finally, as promised, we give our explicit example. 
\begin{theo}
\label{Tinfiniteexample}
Let $$f(z) = \exp(e^z - 1) = \sum_{n=0}^\infty b_n z^n \ \ \ \text{ and } \ \ \ g(z) = \sum_{n=0}^\infty b_{n^2} z^{n^2}.$$ Then $g$ is a {\logreg} {\tef} of infinite order, and $A_R(g)$ is a {\spw}, where $R > 0$ is such that $M(r, g) > r$ for $r \geq R$.
\end{theo}
\proc{Proof.}
We can see that $f$ is {\logreg} because $\phi(t) = \log M(e^t, f) = (e^{e^t} - 1)$ and $$\frac{\phi'(t)}{\phi(t)} > e^t \geq \frac{2}{t}, \qfor t \geq 1.$$

Conditions (\ref{nkdecrease}) and (\ref{gaps}) are 
satisfied, and the coefficients $b_n$ are decreasing because the function $z \mapsto e^z - 1$ satisfies the conditions of Lemma~\ref{decreasinglemma}. Hence, by Theorem~\ref{infinitetheorem},  $g$ is {\logreg} and $A_R(g)$ is a {\spw}.

Finally, $f$ has infinite order. We see from the proof of Theorem~\ref{infinitetheorem} that $f$ and $g$ satisfy (\ref{gandfeq}). Hence, by Lemma~\ref{inforderlemma}, $g$ has infinite order.
\ep\medbreak
Clearly this approach can be used with the function $f$ of Theorem~\ref{Tinfiniteexample} to give a class of functions with $A_R(f)$ {\spws}, by suitably selecting terms from the power series of $f$. We can also use Lemma~\ref{decreasinglemma} to find other {\tef}s which can be manipulated in this way to give further classes of examples.

\proc{Remark.}
We note in passing that, in Theorem~\ref{Tinfiniteexample}, $b_n = B_n/n!$, where $(B_n)$ are the Bell numbers (see, for example, \cite{MR1523147}). Thus, by (\ref{dec1}), we have $$B_{n+1} \leq (n+1)B_n, \qfor n \geq 1.$$ In fact the more precise estimates $$2B_n < B_{n+1} < (n+1) B_n, \qfor n\geq 2,$$ hold (see \cite[Corollary 8]{MR2352650}). These can be deduced in a straightforward way from the identity
\begin{equation}
\label{belleq}
B_{n+1}=\sum_{k=0}^{n}{{\binom{n}{k}}B_k}, \qfor n\geq 0,
\end{equation}
which follows from (\ref{baeq1}).
\medbreak

%
%
\section{Transcendental entire functions with no unbounded Fatou components}
\label{Sbounded}
Baker \cite{MR621564} posed the question of whether the Fatou set of a transcendental entire function of sufficiently small growth can have any unbounded components; see the survey article on this question by Hinkkanen \cite{MR2458805}. By \cite[Theorem 1.5(b)]{fast}, when $A_R(f)$ is a {\spw}, $F(f)$ has no unbounded components. Hence all the examples in this paper have no unbounded Fatou components. In this section we give two results on functions with no unbounded Fatou components, which generalise existing results of this type.

%
%
Our first class of functions with no unbounded Fatou components consists of functions formed by composition of {\psireg} functions.
\begin{theo}
\label{TPbounded}
Let $f_1, f_2, \ldots, f_k$ be {\tef}s which satisfy Lemma~\ref{compositionlemma}(a) for some $m>1$. Suppose that, for all $j \in \{1, 2, \ldots, k\}$, $f_j$ is {\psireg}, with {\psiregfun} $\psi_m$ for $m>1$. Let $g = f_1 \circ f_2 \circ \cdots \circ f_k$. Then every component of $F(g)$ is bounded.
\end{theo}
\proc{Proof.}
By Theorem~\ref{compthpsi}, $A_R(g)$ is a {\spw}, and the result follows by \cite[Theorem 1.5(b)]{fast}.
\ep\medbreak
To compare Theorem~\ref{TPbounded} to previous results, we need the following lemma, part of \cite[Theorem 6]{MR2544754}. This gives a sufficient condition for a {\tef} to be {\psireg}.  We note that although, for other reasons, the full statement of \cite[Theorem 6]{MR2544754} supposes order less than ${\half}$, finite order is sufficient for the proof of this part of the result.
\begin{lemma}
\label{Lpsireg}
Let $f$ be a {\tef} of finite order. Suppose that there exist $n\isnatural$ and $q\in(0,1)$ such that
\begin{equation}
\label{newineq}
M(r,f) \geq \exp^{n+1}((\log^n r)^q), \qfor \text{large } r.
\end{equation}
Then $f$ is {\psireg} with {\psiregfun} given, for all $m>1$, by $$\psi_m(r) = \exp^n((\log r)^p), \quad \text{where } pq > 1.$$
\end{lemma}
The next result now follows from Lemma~\ref{Lpsireg} and Theorem~\ref{TPbounded}.
\begin{cor}
\label{Ccompgen}
Let $f_1, f_2, \ldots, f_k$ be {\tef}s of finite order which satisfy Lemma~\ref{compositionlemma}(a) for some $m>1$. Suppose that there exist $n\isnatural$ and $q\in(0,1)$ such that, for all $j \in \{1, 2, \ldots, k\}$,
\begin{equation}
M(r,f_j) \geq \exp^{n+1}((\log^n r)^q), \qfor \text{large } r.
\end{equation}
Let $g = f_1 \circ f_2 \circ \cdots \circ f_k$. Then every component of $F(g)$ is bounded.
\end{cor}

Rippon and Stallard, in \cite[Theorem 6]{MR2544754}, showed that if $f$ is a {\tef} of order less than $\half$, which satisfies (\ref{newineq}) for some $n\isnatural$ and $q\in(0,1)$, then $f$ has no unbounded Fatou components. By Lemma~\ref{orderhalflemma} this is included in Corollary~\ref{Ccompgen}, with $k=1$.

Corollary~\ref{Ccompgen}, with $n=1$, includes a result of Singh in \cite[Theorem 1]{MR2345508}. (We note that the statement of \cite[Theorem 1]{MR2345508} omits the requirement of finite order, but this was assumed in the proof of \cite[Lemma 1]{MR2345508}.)

%
%
Our second class of functions with no unbounded Fatou components consists of functions formed by composition of {\ef}s, one of which is {\logreg}.
\begin{theo}
\label{Tbounded}
Let $f_1, f_2, \ldots, f_k$ be {\ef}s. Suppose that, for all $j\in\{1,2,\ldots,k\}$, $f_j$ satisfies Lemma~\ref{compositionlemma}(a) for some $m>1$. Suppose also that, for some $j\in\{1, 2, \ldots, k\}$, $f_j$ is a {\logreg} {\tef}. Let $g = f_1 \circ f_2 \circ \cdots \circ f_k$. Then every component of $F(g)$ is bounded.
\end{theo}
\proc{Proof.}
By Theorem~\ref{compth}, $A_R(g)$ is a {\spw}, and the result follows by \cite[Theorem 1.5(b)]{fast}.
\ep\medbreak
As noted in Section~\ref{Sprelim}, this result differs from Theorem~\ref{TPbounded} in that only one function in the composition needs to satisfy the regularity condition and be transcendental.

The final result follows from Theorem~\ref{Tbounded} and Lemma~\ref{orderhalflemma}.
\begin{cor}
\label{Cbounded}
Let $f_1, f_2, \ldots, f_k$ be {\tef}s of order less than $\half$. Suppose that, for some $j \in \{1, 2, \ldots, k\}$, $f_j$ is {\logreg}. Let $g~=~f_1 \circ f_2 \circ \cdots \circ f_k$. Then every component of $F(g)$ is bounded.
\end{cor}
This corollary generalises a result of Anderson and Hinkkanen in \cite[Theorem 2]{MR1452790}, which states that if a {\logreg} function has order less than ${\half}$, then it has no unbounded Fatou components. Anderson and Hinkkanen's result is included in Corollary \ref{Cbounded} with $k=1$.

Cao and Wang \cite{MR2105780} developed a similar result to Corollary~\ref{Cbounded}, concerning composition of {\tef}s. They showed that if $g~=~f_1~\circ~f_2~\circ~\cdots~\circ~f_k$, where $f_1, f_2, \ldots, f_k$ are {\tef}s of order less than ${\half}$, at least one of which has positive lower order, then $g$ has no unbounded Fatou components. By Lemma~\ref{logreglemma}, Cao and Wang's result is included in Corollary~\ref{Cbounded}. We note that it is possible to construct a class of {\logreg} functions of lower order zero and any given finite order, in particular order less than $\half$. This shows that there are situations in which Corollary \ref{Cbounded} can be applied but not the result of \cite{MR2105780}.

\ack The author is grateful to Phil Rippon and Gwyneth Stallard for their patient, invaluable help with this paper.
%
%
\bibliographystyle{acm}
\bibliography{Main}
\end{document}